\documentstyle[11pt,leqno,here]{article}
\input amssym.def
\input amssym.tex
\input psfig

\newcommand{\beql}[1]{\begin{equation}\label{#1}}
\newcommand{\eeq}{\end{equation}}

\makeatletter
\@addtoreset{figure}{section}
\def\thefigure{\thesection.\@arabic\c@figure}
\@addtoreset{table}{section}
\def\thetable{\thesection.\@arabic\c@table}

\def\@sect#1#2#3#4#5#6[#7]#8{\ifnum #2>\c@secnumdepth
     \def\@svsec{}\else
     \refstepcounter{#1}\edef\@svsec{\csname the#1\endcsname.\hskip .75em }\fi
     \@tempskipa #5\relax
      \ifdim \@tempskipa>\z@
        \begingroup #6\relax
          \@hangfrom{\hskip #3\relax\@svsec}{\interlinepenalty \@M #8\par}%
        \endgroup
       \csname #1mark\endcsname{#7}\addcontentsline
         {toc}{#1}{\ifnum #2>\c@secnumdepth \else
                      \protect\numberline{\csname the#1\endcsname}\fi
                    #7}\else
        \def\@svsechd{#6\hskip #3\@svsec #8\csname #1mark\endcsname
                      {#7}\addcontentsline
                           {toc}{#1}{\ifnum #2>\c@secnumdepth \else
                             \protect\numberline{\csname the#1\endcsname}\fi
                       #7}}\fi
     \@xsect{#5}}
\def\@begintheorem#1#2{\it \trivlist \item[\hskip \labelsep{\bf #1\ #2.}]}
\def\section{\@startsection {section}{1}{\z@}{-3.5ex plus -1ex minus
 -.2ex}{2.3ex plus .2ex}{\normalsize\bf}}

\makeatother

\begin{document}
\begin{center}
{\Large {\bf Repeated Patterns of Dense Packings 
\\of Equal Disks in a Square \\}} 
\vspace{1\baselineskip}
{\em R. L. Graham},~~ rlg@research.att.com \\
{\em B. D. Lubachevsky}, ~~bdl@research.att.com \\
\vspace*{1\baselineskip}
AT\&T Bell Laboratories  \\
Murray Hill, New Jersey 07974, USA \\
\vspace{1\baselineskip}
%Submitted: January 17, 1996; Accepted: April 20, 1996 \\
\vspace{1.5\baselineskip}
{\bf ABSTRACT}
\vspace{.5\baselineskip}
\end{center}

\setlength{\baselineskip}{0.995\baselineskip}

We examine sequences of dense packings
of $n$ congruent non-overlapping disks inside a square
which follow specific patterns as $n$ increases along
certain values, $n = n(1), n(2),... n(k),...$.
Extending and improving previous work
of Nurmela and  \"Osterg\aa rd \cite{NO} where
previous patterns for $n = n(k)$ of the form 
$ k^2$, $ k^2-1$, $k^2-3$, $k(k+1)$, and $4k^2+k$ were observed,
we identify new patterns for $n = k^2-2$
and $n = k^2+ \lfloor k/2 \rfloor$.
We also find denser packings than those in \cite{NO}
for $n =$21, 28, 34, 40, 43, 44, 45, and 47.
In addition, we produce what we conjecture to be
optimal packings for
$n =$51, 52, 54, 55, 56, 60, and 61.
Finally, for each identified sequence $n(1), n(2),... n(k),...$
which corresponds to some specific repeated pattern,
we identify a threshold index $k_0$, for which
the packing appears to be optimal for $k \le k_0$,
but for which the packing is not optimal (or does
not exist) for $k > k_0$.

\setlength{\baselineskip}{1.2\baselineskip}
\section{Introduction}
\begin{figure}[htb]
\centerline{\psfig{file=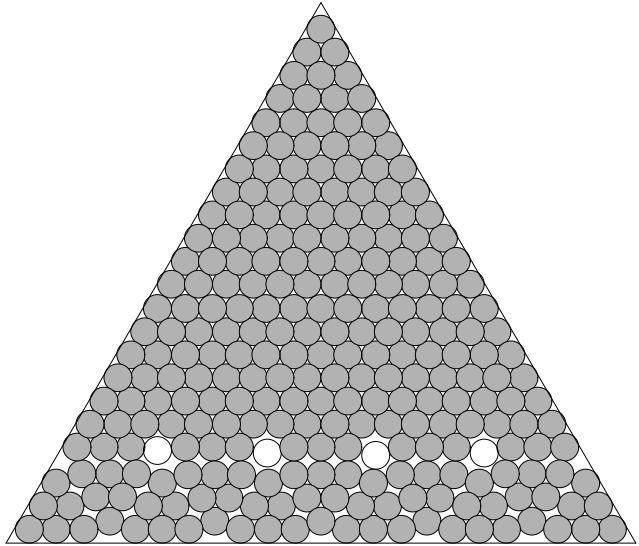,width=3.5in}}

\caption{The conjectured densest packing of $n = 256$ disks inside an equilateral triangle,
a member of the series $n = n_p (k) = ~\Delta((k+1)p-1) + (2p+1) \Delta(k)~$
for $p = 5$ and $k = 3$, where $\Delta(m)~=~m(m+1)/2$. 
The $n-p+1 = 252$ shaded disks can not move (they are ``solid''),
the $p-1 = 4$ non-shaded disks are free to move within their local confines
(they are ``rattlers'').
The densest packings of $n$ disks for all
checked values of the form $n = n_p (k), p = 1,2,..., k = 1,2,...$,
have this pattern consisting of one triangle of side $(k+1)p-1$
and $2p+1$ triangles of side $k$ with $p-1$ rattlers 
that are ``falling off'' the larger triangle.}
\label{t256}
\end{figure}

\hspace*{\parindent}
In a previous paper \cite{GL1}, the authors observed
the unexpected occurrence of repeating ``patterns"
of dense (and presumably optimal) packings of
$n$ equal non-overlapping disks inside an equilateral triangle
(see Fig.~\ref{t256} for and example). 
It is natural to investigate this phenomenon for other
boundary shapes. In particular, this was done by the
authors \cite{LG1} for the case of $n$ disks in a circle.
However, in contrast to the case of the equilateral
triangle where the patterns appear to persist for
arbitrarily large values of $n$, for the circle
the identified packing patterns cease to be optimal as
the number of disks exceeds a certain threshold.

In this note we describe the situation for the square.
In a recent paper, Nurmela and  \"Osterg\aa rd \cite{NO} present
various conjectured optimal packings of
$n$ equal disks in a square for up to 50 disks. They also point
out certain patterns that occur there. 
By using a packing procedure different 
from theirs, we improve on their best packings for
$n=21$, 28, 34, 40, 43, 44, 45, and 47.
We conjecture that these new packings are
optimal, as are the new packings we give for
$n=51$, 52, 54, 55, 56, 60, and 61.

We confirm all the repeated patterns
mentioned in \cite{NO},
specifically, for
$n = k^2$, $ k^2-1$, $k^2-3$, and $k(k+1)$,
and we identify two new patterns, namely, for $n = k^2-2$
and $n = k^2+ \lfloor k/2 \rfloor$
The latter pattern incorporates the
packings of $n = 4k^2+ k$ disks 
as identified in \cite{NO}.

It was found in \cite{NO} that the obvious ``square''
pattern of packings of $n = k^2$ disks
becomes non-optimal for $n > n_0 = 36$.
This was done by presenting a configuration of $k^2=49$ disks
with the diameters larger than $m = 1/(k-1)$.
The latter $m$ is disk diameter in the packing that obeys the pattern.
The standard unit of measure
used in most papers on the subject is
the side of the smallest square that contains the {\em centers}
of disks.
We repeat the same procedure for the other patterns identified both
in \cite{NO} and in the present paper.
Namely, 
for each pattern we state the rule of its formation
which allows us to compute the corresponding value of $m = m(n)$.
Then we pinpoint
the $n_0$ that belongs to the series and such that the packing of
$n_0$ disks constructed according to the rule is
(presumably) optimal, but for which $m(n_1)$ for the
next value $n_1 > n_0$ in the series when the packing is also 
constructed according to the rule
is {\em worse} than a certain {\em challenger} disk configuration
(which may or may not be a solid packing).

In this manner, we confirm observation in \cite{NO}
that the best packings of $n = k^2-1$ disks loose their characteristic
pattern for $n \ge n_1 = 48$ disks.
We also found that
although it was not stated in \cite{NO},
the packing of 47 disks presented there 
(as well as our better packing of 47 disks)
challenges the series $k^2-2$.
Thus, 
the pattern of the series $n = k^2-2$ 
becomes non-optimal for $n > n_0 = 34$,
and that of the series $n = k^2-1$ for $n > n_0 = 35$.
We also found challenger disk configurations or packings 
for other patterns for values of $n > 50$ which were
not identified in \cite{NO}. 
Namely:
$n_0 = 56$ for the series $k(k+1)$ 
(with the challenger $n_1 = 72$),
and $n_0 = 61$ for the series $k^2 -3$
(with the challenger $n_1 = 78$).
The situation for the series $n = n(k) = k^2+ \lfloor k/2 \rfloor$
($n = 5$, 10, 18, 27, 39, 52, 68,...) is more complex:
the pattern exists for $5 \le n(k) \le 52$ as a solid packing
and is probably optimal 
for these $n(k)$ except the case $n=n(3)=10$
which is the subject of several publications 
(\cite{G}, \cite{Sch}, \cite{Schl}, \cite{Val}).
For $n = n(8) = 68$ the configuration constructed according
to the pattern rule has a slight disk overlap,
i.e., it {\em does not exist} as a disk packing,
and the overlap persists for all $n = n(k) > 68$.

\begin{figure}[htb]
\centerline{\psfig{file=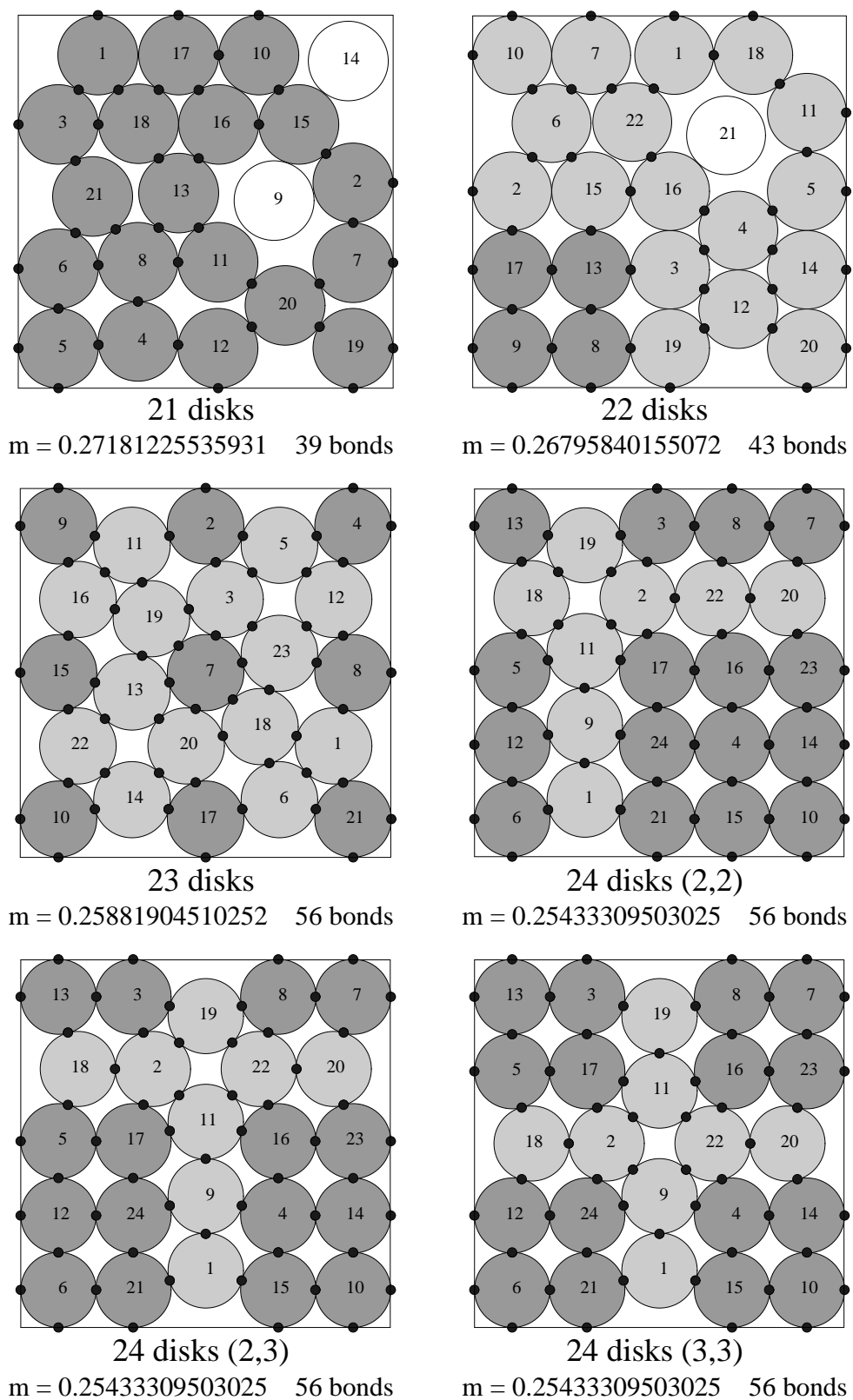,width=5.3in}}
\caption{The densest packings found of 21 to 24
disks}
\label{t21-24}
\end{figure}

In the geometric packing problem, 
progress in  {\em proving} lags that of {\em conjecturing}.
Thus, we should warn the reader that almost all 
our statements are conjectures; they are based on
computer experimentations with the so-called
``billiards'' simulation algorithm \cite{L}, \cite{LS}.
In all series except $k^2-3$, the construction rule
we found for generating a pattern for a given $n$
is a finite procedure and $m(n)$ can be
expressed as the root of a well-defined polynomial. The existence
of the packing of a given pattern for a fixed $n$,
even if not the pattern's optimality
(when appropriate), should be considered proven. 
However, for the series $k^2-3$
we have to resort on an infinite simulation procedure
\cite{L}, \cite{LS} to construct a pattern for an arbitrary $n$
and to compute $m(n)$. Hence even the pattern's existence
as a solid packing is a conjecture here.
We point out that some of the proposed methods
in the literature attempt to prove a packing optimal
or, at the least, prove that a packing with particular parameters
exists. Usually to fulfill this task,
the packing must be actually presented, even if only as
a conjecture.
Thus, we try to present our conjectured packings
in a verifiable and reproducible form;
we provide 14 decimal digits of accuracy for its parameter $m$ 
and clearly identify the connectivity pattern
(touching disk-disk and disk-wall pairs).
Some previous papers provided disk coordinates in the presented packings.
\footnote{
We were unable to reproduce the best
packing of 21 disks for which \cite{MFP}
provides the diameter $m = 0.27181675$
but no other data, e.g. no contact diagram.
Our best packing of 21 disks has a smaller diameter
(see Fig. \ref{t21-24}).}
The interested reader can contact the authors directly
for the coordinates (since that would otherwise take up 
too much space in the paper).

\section{Packings}
\begin{figure}[htb]
\centerline{\psfig{file=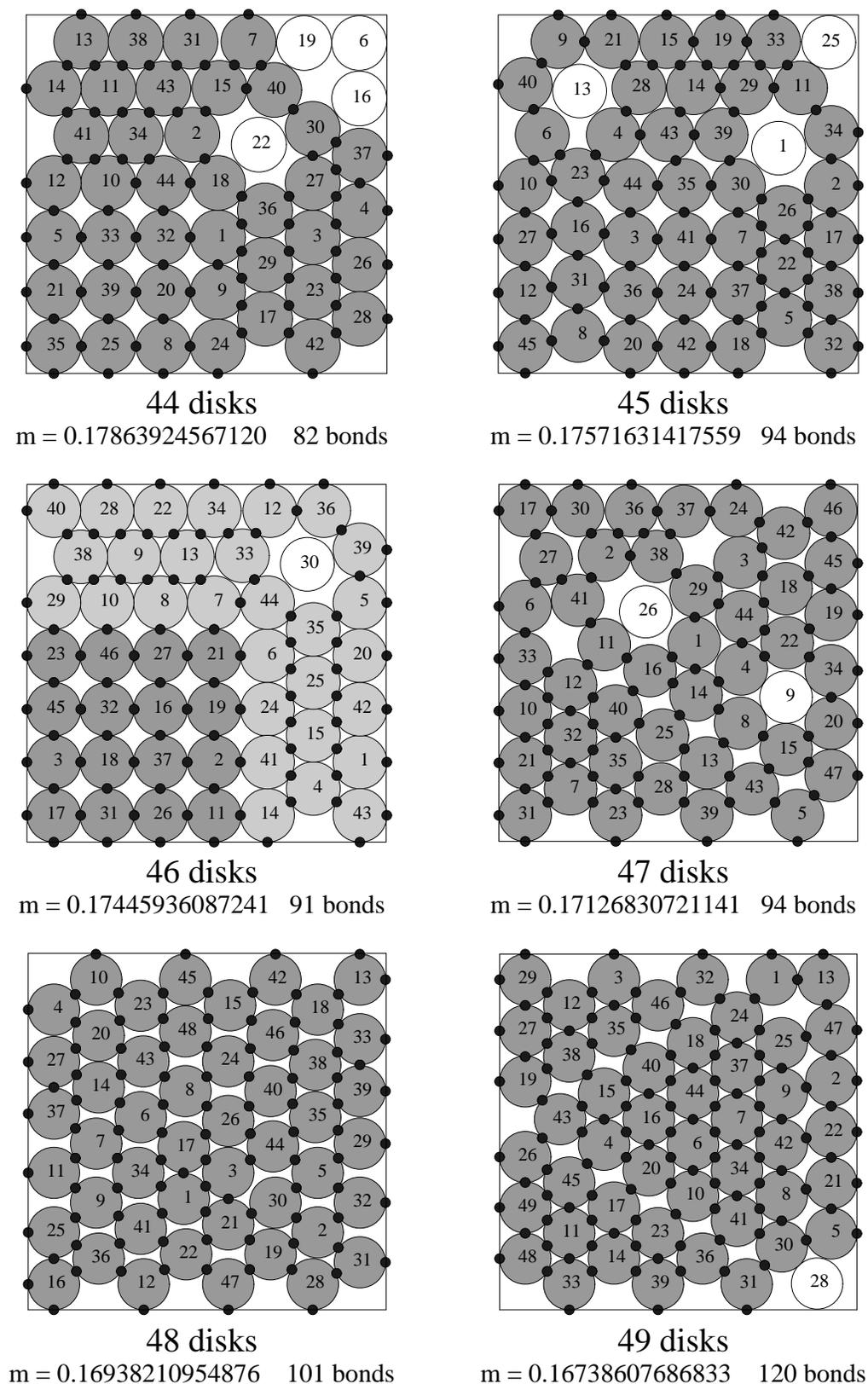,width=5.5in}}
\caption{The densest packings found of 44 to 49
disks}
\label{t44-49}
\end{figure}

\hspace*{\parindent}
The parameter $m$ supplied with each packing is the
ratio of the disk diameter to the side length of the
smallest square that contains the disk centers.
{\em Bonds} or contact points mark disk-disk or disk-boundary contact.
In the packing diagrams, bonds are indicated by black dots.
Most of the packings presented are conjectures.
This means that a proof is required 
not only for their optimality, if any,
but even for their existence.
Thus, a bond implies a conjecture that the corresponding
distance is zero while the absence of a bond implies
a conjecture
that the distance is strictly positive.
We placed or did not place
a bond between two disks or a disk and a boundary
based on the numerical evidence:
the bond was placed
when the corresponding distance was less than $10^{-12}$ of the disk
diameter.
Such a choice of a threshold is supported by the existence
of a well-formed gap between a bond and a no-bond situation:
In all cases when the bond is not present between apparently
touching surfaces, the computed distance is at least $10^{-7}$
of the disk diameter,
and, except for the packing of 47 disks in Fig.\ref{t44-49},
it is at least $10^{-5}$.
The existence of this gap also testifies to that in all the packings
the double precision resolution we employed for the computations 
sufficed.
All {\em solid} disks, i.e.,  those that can not move, 
are shaded in the packing diagrams; 
the non-shaded disks are {\em rattlers}---
they are free to move within their confines.
Different shadings of disks in some packings
and a unique numerical label for each disk on a diagram are
provided to facilitate the discussion,
These are not part of the packings.

\begin{figure}[htb]
\centerline{\psfig{file=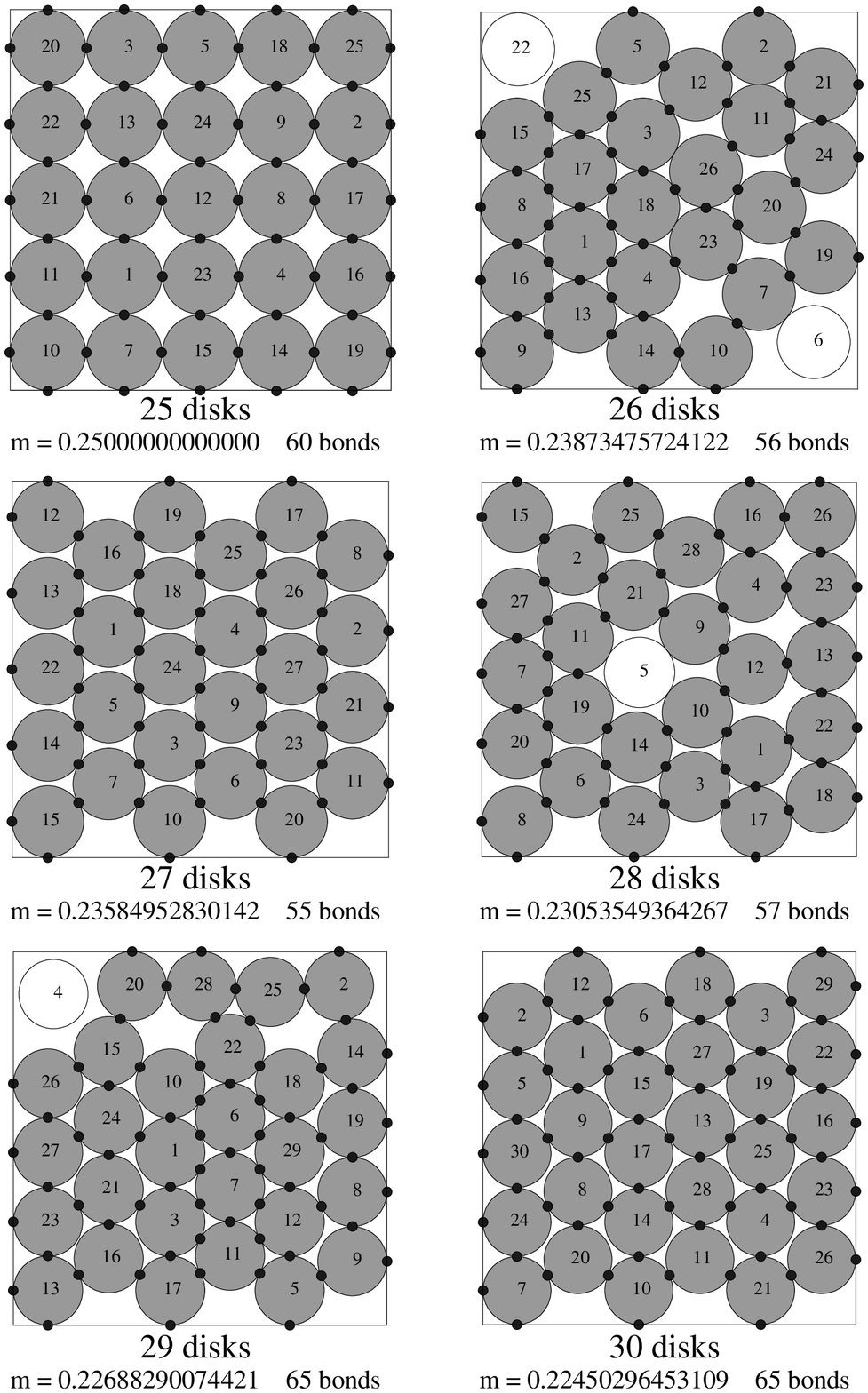,width=5.5in}}
\caption{The densest packings found of 25 to 30 
disks}
\label{t25-30}
\end{figure}

\begin{figure}[htb]
\centerline{\psfig{file=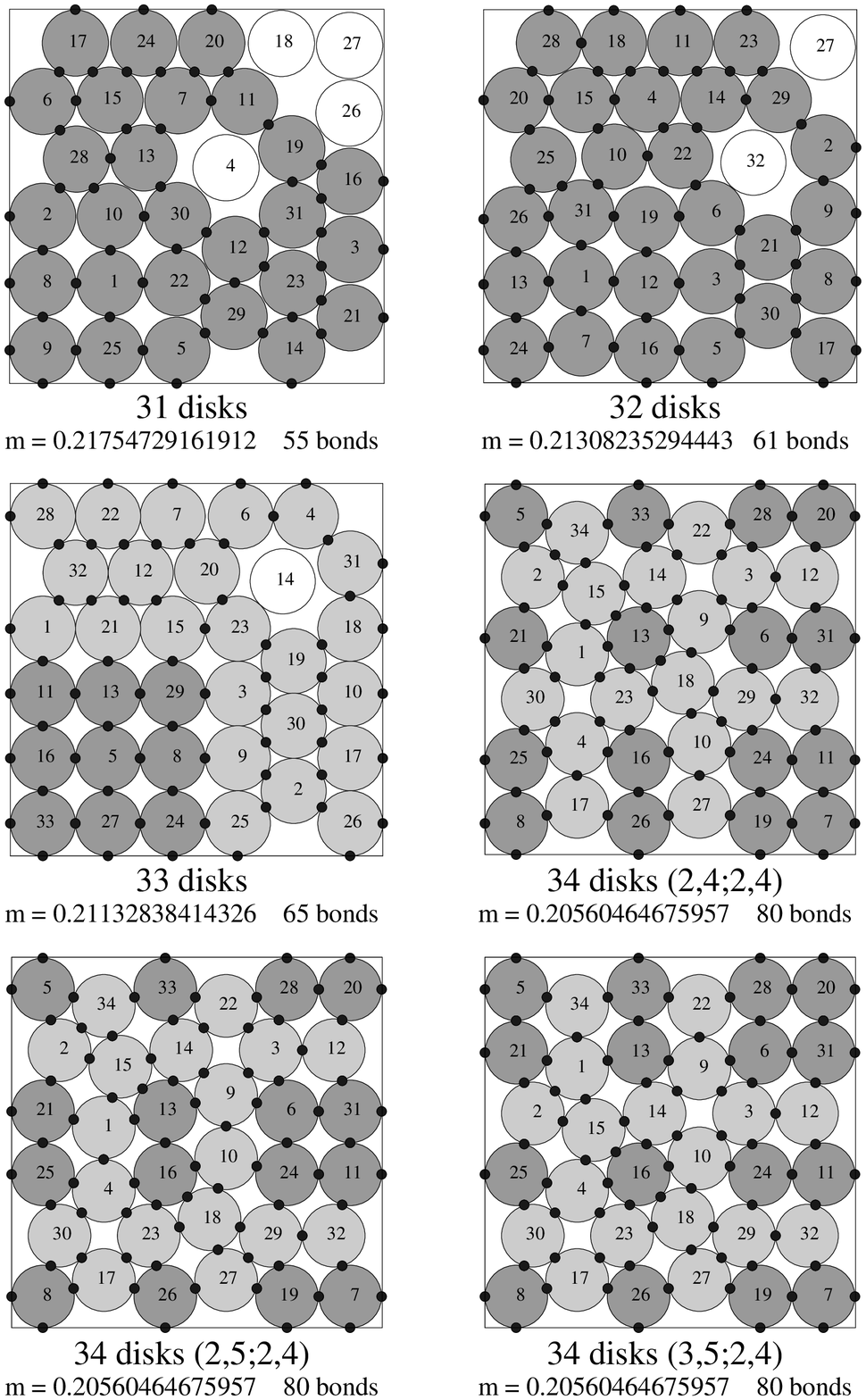,width=5.5in}}
\caption{The densest packings found of 31 to 34
disks}
\label{t31-34}
\end{figure}

\begin{figure}[htb]
\centerline{\psfig{file=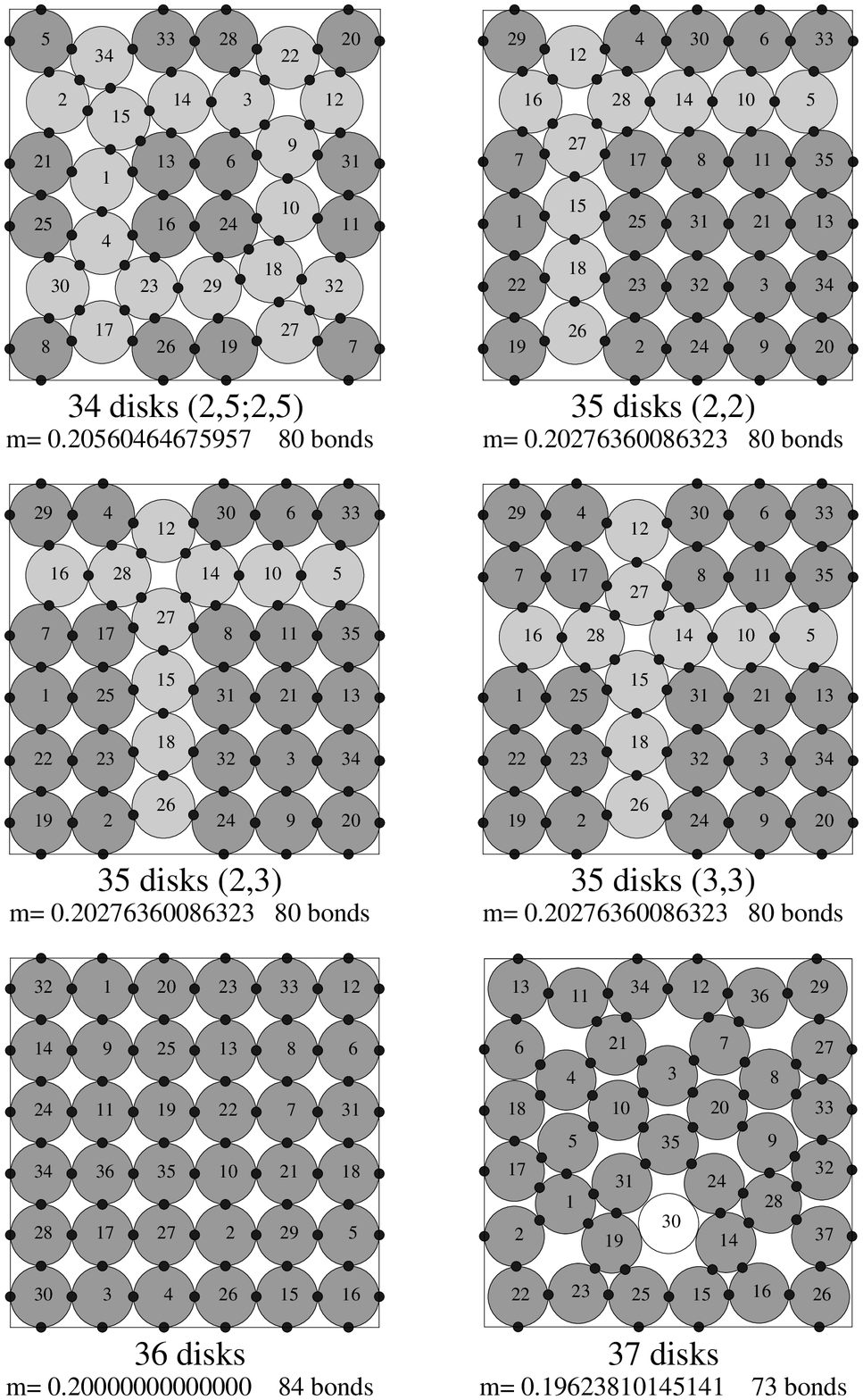,width=5.5in}}
\caption{The densest packings found of 34 to 37
disks}
\label{t34-37}
\end{figure}

\begin{figure}[htb]
\centerline{\psfig{file=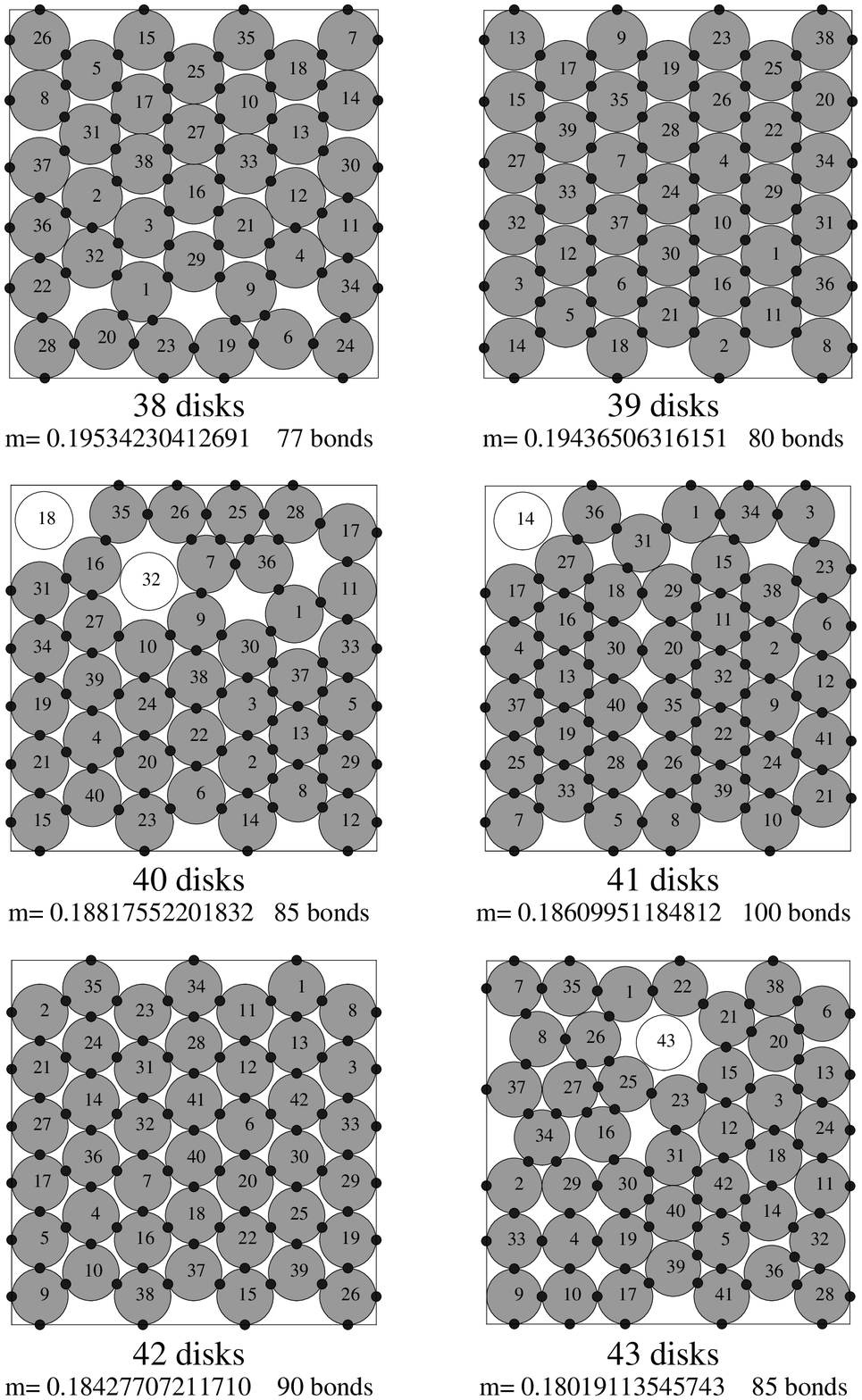,width=5.5in}}
\caption{The densest packings found of 38 to 43
disks}
\label{t38-43}
\end{figure}

\begin{figure}[htb]
\centerline{\psfig{file=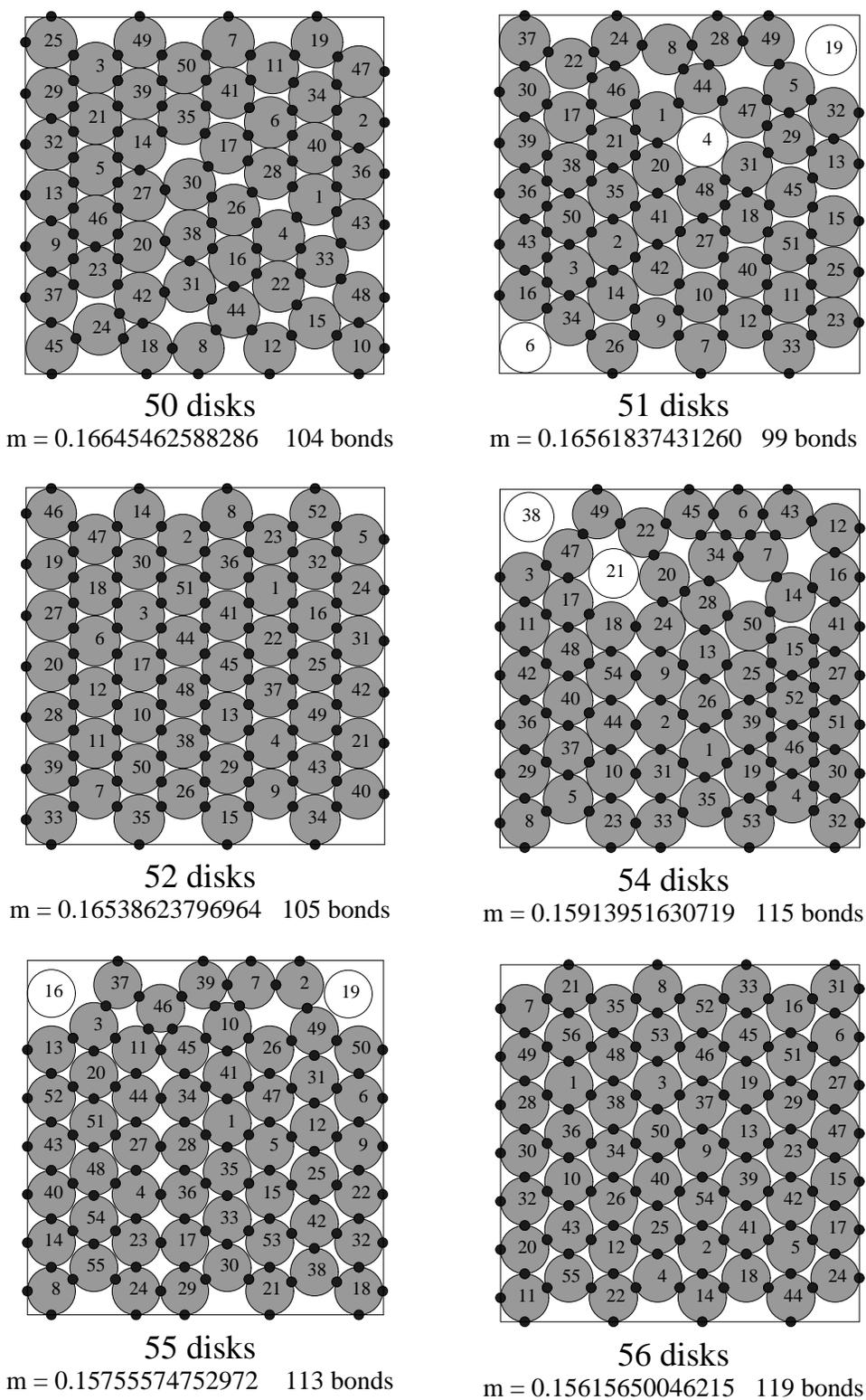,width=5.2in}}
\caption{The densest packings found of 50 to 52, and 54 to 56 disks}
\label{t50-56}
\end{figure}

\begin{figure}[htb]
\centerline{\psfig{file=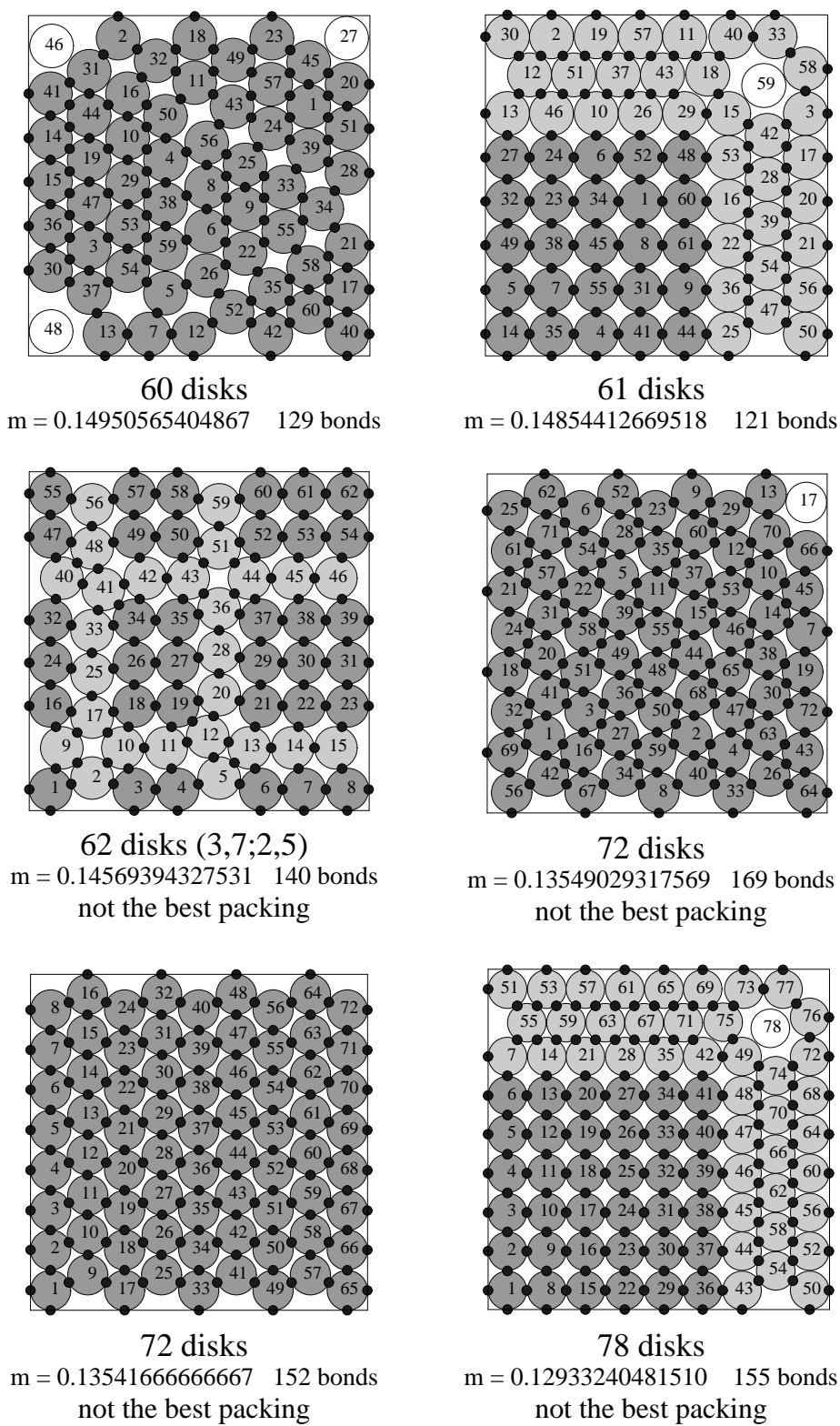,width=5.0in}}
\caption{The densest packings found of 60 and 61 disks and 
inferior packings of 62, 72, and 78 disks}
\label{t60-x}
\end{figure}

\begin{figure}[htb]
\centerline{\psfig{file=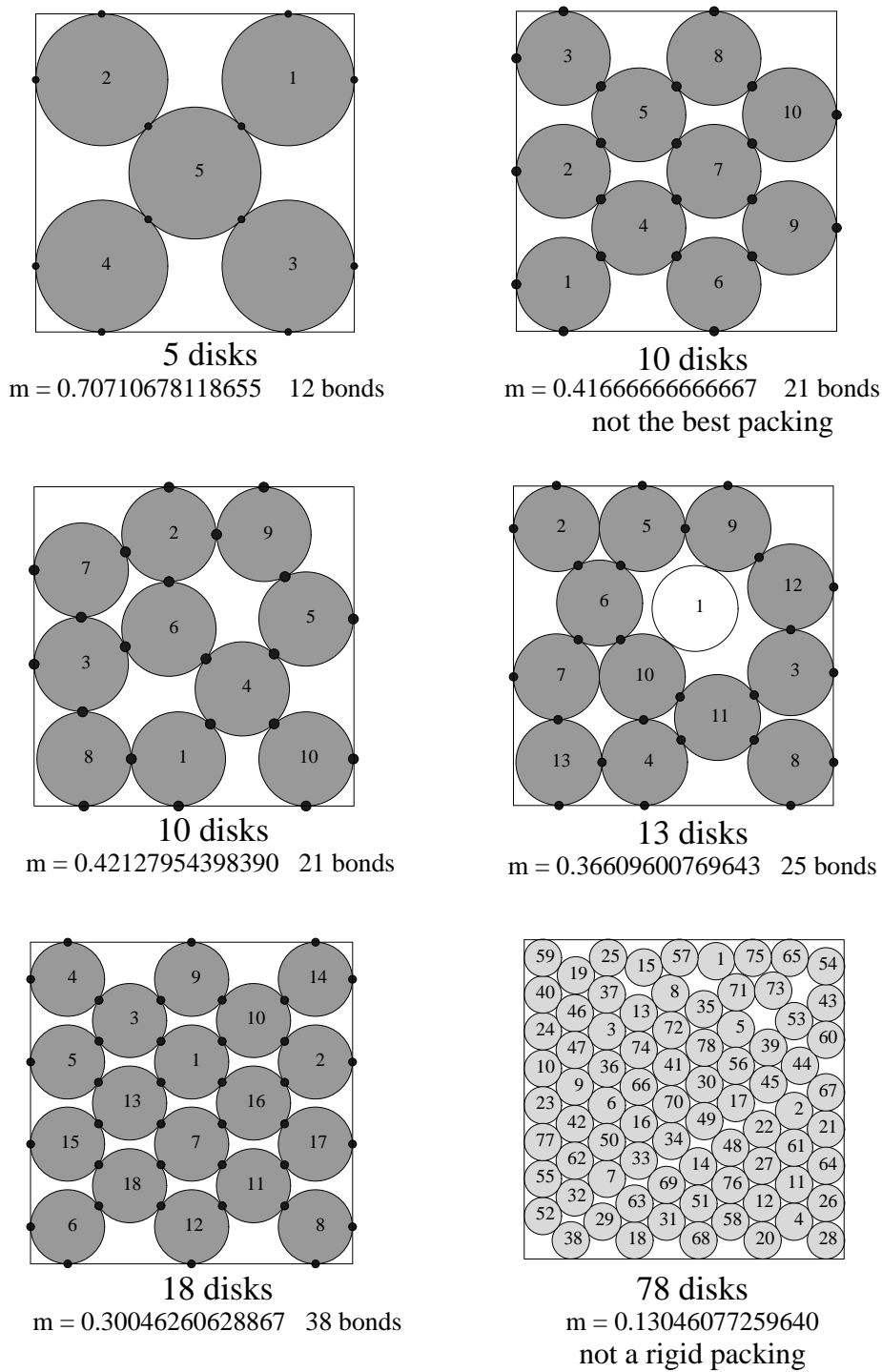,width=5.0in}}

\caption{The densest packings of 5, 10, 13, and 18 disks,
an inferior packing of 10 disks, and a challenger configuration of 78 disks}
\label{t5-x}
\end{figure}

\clearpage

\section{Series of packings with similar patterns}
%\hspace*{\parindent}
\paragraph{${\bf n = k^2}$}.
The packings of this series have an obvious
square pattern with disk diameter $m = 1/(k-1)$.
The pattern is optimal for $n = 1, 4, 9$, and 16, 
and probably $25$ and $36$.
For $k_1 = 7$ this yields value of $m = 1/6$ which is smaller
than the experimental diameter of the challenger packing
of 49 disks in Fig. \ref{t44-49}.
\paragraph{${\bf n = k^2-1}$}. The pattern can be viewed as a square
of $(k-1)^2$ disks arranged in straight
rows and columns into which one 
row and one column of ``shifted'' disks are inserted.
A ``shifted'' row or column each has $k-1$ disks.
In the diagrams for $k = 5$ ($n = 24$, Fig. \ref{t21-24})
and $k = 6$ ($n = 35$, Fig. \ref{t34-37}), the
``straight'' disks are shaded more heavily
than the ``shifted'' ones.
Several equivalent packings are obtained
by inserting the shifted row and column in different places
among the straight ones.
A pair of insertion position $(i,j), 2 \le i,j \le k-1$,
identifies the packing, e.g., in packing of 24 disks (2,3)
(Fig. \ref{t21-24}),
the shifted row is the second and the shifted column is the third
(counting from the top left corner).
Among the $(k-2)^2$ packings thus obtained,
there are many congruent pairs.
The pattern is fully developed at $n = 8$ (one packing)
and remains optimal
for $n = 15$ (one packing),
$n=24$ (three equivalent packings), 
and $n = n_0 = 35$ (three equivalent packings).
For $n_1 = 48$ the pattern looses its optimality,
which can be shown as follows.
The angle at the center of disk 13 in packing (2,2) of 24 disks
(see Fig. \ref{t21-24} ) in the triangle formed by centers of disks 
13, 19 and 3 is $15^o$ with $\cos 15^o = .5 \sqrt{ 2 + \sqrt{3} } $.
Similar angles are the same in all the packings
of the series which easily implies
$m = 1/(k - 3 + \sqrt {2 + \sqrt {3}} )$.
For $k_1 = 7$ this yields value of $m = 0.168581424$ which is smaller
than the experimental diameter of the challenger packing
of 48 disks in Fig. \ref{t44-49}.
\paragraph{${\bf n = k^2-2}$}. This pattern 
is similar to the pattern of series $k^2-1$,
only here there are two shifted rows and two shifted columns.
The (non-optimal) packing of 62 disks depicted in Fig. \ref{t60-x}, 
shows this pattern for $k = 8$.
There are $(k-3)(k-4)/2$ possible ways to insert a pair of shifted
rows $(i_1 ,i_2)$ 
because of the restrictions $1 < i_1 < i_2 < k$
and similarly for the columns $(j_1 ,j_2)$.
Hence, there are
$((k-3)(k-4)/2)^2$ possible different index sets ($i_1 ,i_2 ;j_1 ,j_2)$
for each of which we can construct an equivalent packing.
Many of these are congruent.
The pattern is fully developed at $k = 5$ ($n = 23$, one packing) 
and remains optimal
for only one more value $k = k_0 = 6$ 
($n_0 = 34$, four equivalent packings.).
Here $m = 1/(k - 5 + 2 \sqrt {2 + \sqrt {3}})$. For $k = k_1 = 7$
this yields value $m = 0.1705406887$ which is smaller than that of
the challenger packing of $n_1 = 47$ disks in Fig. \ref{t44-49}.
\paragraph{${\bf n = k^2-3}$}. This pattern is represented
by conjectured dense packings of 22, 33, 46, and 61 disks 
(Figs. \ref{t21-24},  \ref{t31-34}, \ref{t44-49},  \ref{t60-x}).
A feature of the pattern is a (shaded more heavily on the pictures)
densely packed ``straight'' square of $(k-3)^3$ disks
in the bottom left corner.
It follows that the pattern is not optimal for $k \ge 10$ ($n \ge 97$)
because the ``straight'' square itself isn't.
However, even for a smaller
$n = n_1 = 78$ we were able to get a disk configuration
that challenges the pattern as presented in Fig. \ref{t60-x},
although this configuration is a not fully
formed solid disk packing (see Fig. \ref{t5-x}).

Looking at the (non-optimal) packing of 78 disks (see Fig. \ref{t60-x}),
the pattern can be further described as 
three alternating columns at the right
and three alternating rows at the top
with one rattler at the top right corner.
In each of the three additional rows 
and
the three additional columns, most disks
are not touching each other.
The exceptions are 
pairs 42-49 in the bottom additional row,
73-77 in the top additional row
and 76-72 in the right additional column
(Fig. \ref{t60-x}).
Also there is almost full contact for pairs of disks
between adjacent additional rows and columns,
except for the pair 42-75 which are not touching
between the first and second additional
rows.
These features are identical 
in all the packings of the pattern.

Each such packing can be obtained by ``tightening''
a certain 
configuration $C$ described
as follows:
In $C$, disks 77 and 76 are not touching
each other but
all disks
in additional rows and columns are touching their neighbors.
($C$ is not a solid packing.)
%In $C$ we have $m = 1/(k - 4 + \sqrt {3})$.
%This $m$ 
%can serve as a lower bound for the densest packing
%of $k^2-4$ disks.
We
take $C$ as an initial condition
for the ``billiards'' packing algorithm \cite{L}, \cite{LS}.
To ``tighten'' the configuration,
the ``billiards'' algorithm
allows the disks to move chaotically in the square
without overlaps
while their diameter
increases
at a common rate
until no further growth is possible.
It is remarkable, that for each $n=k^2-3$,  $k = 5,6,7,8$ and 9,
this chaotic negotiation
always converges (independently of the initial velocities)
to the described pattern with parameter $m$ identical to double
precision
in all runs involving the same $n$. 
Moreover, when instead of the configuration $C$,
we start with a random initial configuration 
and zero initial diameters of the disks, the same configuration
with the same precision results in the runs that achieve the largest $m$,
- but only for $k = 5,6,7$ and 8, excluding value $k=9$.
Also, for $k = 4$ ($n = 13$) which seems to be the smallest $n$
for which the pattern may exist,
both methods,
the one beginning with a zero-diameter random configuration
and the one beginning with the configuration $C$, lead to the same packing
which, however, deviates somewhat from the described pattern,
see Fig. \ref{t5-x}
\paragraph{${\bf n = k(k+1)}$}. The pattern
consists of $k+1$ alternating columns with $k$ disks each,
see, e.g., the packing of 30 disks in Fig. \ref{t25-30}.
Following the path of disks 2, 12, 6, 18, 3, 29 in this packing,
we come to the expression $k \cos \alpha$ for the horizontal side of
the square that contains the disk centers,
where the angle $\alpha = \alpha (k)$ is that at the vertex 2
in the triangle 12, 2, 6.
The vertical side of the square is $k-1+ \sin \alpha$, where
the latter expression results from the path 12, 2, 5, 30, 24, 7.
Equating these two expressions for the side we 
determine from the resulting equation:
$\cos \alpha (k) = (k^2-k+\sqrt{2k})/(k^2+1)$ and 
$m = (k \cos \alpha (k) )^{-1}$.

Looking again at the packing of 30 disks in Fig. \ref{t25-30},
circles 1 and 2 are separated by a non-negative distance 
because $\alpha (k) \le 30^o$
or $\cos \alpha (k) \ge \sqrt{3}/2$
for $k = 5$.
The latter inequality holds 
and the pattern exists as a solid packing of non overlapping disks
only for $k \ge 4$, i.e., for $n = 20$, 30, 42, 56, 72,....

We show that the packing of the pattern is not optimal
for $n \ge n_1 = 72$ as follows.
We present another regular pattern of packings of
$k(k+1)$ disks; for $k = 8$
this alternative pattern is depicted in 
the second row and second column in
Fig. \ref{t60-x}
(with one rattler).
The alternative pattern
exists for all values of $n$ for which
the main pattern exists. 
The diameter $\bar{m}=\bar{m}(k)$ of a disk in the alternative pattern
is given by 
$\bar{m}(k) = 1/((k+1/2)\cos \beta (k)+( \sqrt{3}/2) \sin \beta (k))$
where (in the example in Fig. \ref{t60-x})
$\beta (k)$ for $k = 8$ is the angle at disk 56
in the triangle 67, 56, 42.
As before, we 
equate the horizontal and vertical sides of the square:
$k \cos \beta + \cos ( \beta + \pi /3 ) = 
(k-1) \sin ( \beta + \pi /3 ) + \sin \beta$.
From this equation we easily determine $\beta = \beta (k)$ and then
we find that $\bar{m} > m$ for all $n \ge n_1 = 72$
(but $\bar{m} < m$ for $n = 20$, 30, 42, and 56).
Note that the alternative packing,
although it is better than the packing of the pattern for
$n = 72$ disks, is not optimal.
For example, for 72 disks we found experimentally
a disk configuration with an irregular structure (not shown)
that is better than both the main or the alternative pattern.
\paragraph{${\bf n = k^2+ \lfloor k/2 \rfloor }$}. The pattern,
as exemplified by cases $k = 5$ ($n = 27$, Fig. \ref{t25-30})
and $k = 6$ ($n = 39$,  Fig. \ref{t38-43}), consists
of $k+1$ alternating columns, odd columns having $k$ disks each
and even columns having $k-1$ disks each.
As before, with $\alpha$ denoting the angle,
say, at disk 15 in the triangle 7, 15, 10
in the packing of 27 disks in Fig. \ref{t25-30},
we compute the side-length of the square obtained in two different ways:
$k \cos \alpha = 2(k-1) \sin \alpha$.
From this equation we easily determine $\cos \alpha (k)$ and then
$m = 1/(k \cos \alpha (k) )$.
The pattern exists when, as in the example of 27 disks
(Fig. \ref{t25-30}), 
disks 12 and 13 are separated by a non-negative distance.
This occurs when $\sin \alpha (k) \ge 1/2$.
Thus, the pattern exists
only for $k = 2$, 3, 4, 5, 6, and 7 
($n = 5$, 10, 18, 27, 39, and 52)
and those packings (shown in Figs \ref{t5-x}, \ref{t25-30}, \ref{t38-43}, and
\ref{t60-x}) are, indeed, optimal, except the case of $n = 10$: 
proved \cite {GMPW} for $n = 5$, 10, 18
and conjectured for $n = 27$, 39, and 52.

For $k \ge 8$ the ``ideal'' pattern yields overlaps in pairs of disks.
The overlap increases with $k$. 
For $k = 8$ ($n = 68$) the overlap is less than 1\% of the
disk diameter so it is naturally to expect that the optimal
packing will be a small deviation from the main pattern.
Experiments, indeed show that best configurations
are of this sort. 
Unfortunately, none of them is a solid packing,
because
it becomes very difficult to find exactly in which pairs the disks
are touching, and not merely just very close to each other.
%\section{Discussion}
\hspace*{\parindent}

\end{document}